\documentclass[reqno,12pt,twoside]{article}
\usepackage{fullpage}
\usepackage{amssymb}
\usepackage{amsmath}
\usepackage{amsthm}
\newtheorem{theorem}{Theorem}
\newtheorem{problem}{Problem}

\newtheorem{lemma}{Lemma}
\numberwithin{equation}{section}

\begin{document}

\title{\Large On asymptotic bases and minimal asymptotic bases
}\author{\large Min Tang\thanks{Corresponding author. This work was supported by the National Natural Science Foundation of China(Grant
No. 11471017).} and Deng-Rong Ling}
\date{} \maketitle
 \vskip -3cm
\begin{center}
\vskip -1cm { \small
\begin{center}
 School of Mathematics and Statistics, Anhui Normal
University
\end{center}
\begin{center}
Wuhu 241002, PR China
\end{center}}
\end{center}

 {\bf Abstract:} Let $\mathbb{N}$ denote the set of all nonnegative integers and $A$ be a subset of $\mathbb{N}$. Let $h\geq2$ and let
$r_h(A,n)=\sharp \{ (a_1,\ldots,a_h)\in A^{h}: a_1+\cdots+a_h=n\}.$ The set $A$ is called an asymptotic basis of order $h$ if $r_h(A,n)\geq 1$ for all sufficiently large integers $n$. An asymptotic basis $A$ of order $h$ is minimal if no proper subset of $A$ is an asymptotic basis of order $h$. Recently, Chen and Tang resoved a problem of Nathanson on minimal asymptotic bases of order $h$. In this paper, we generalized this result to $g$-adic representations.


\vskip8mm

\section{Introduction}

Let $\mathbb{N}$ denote the set of all nonnegative integers and $A$ be a subset of $\mathbb{N}$. Let $h\geq2$ and let
$$r_h(A,n)=\sharp \{ (a_1,\ldots,a_h)\in A^{h}: a_1+\cdots+a_h=n\}. $$
Let $W$ be a nonempty subset of $\mathbb{N}$. Denote by $\mathcal{F}^{\ast}(W)$ the set of all finite, nonempty subsets of $W$. For any integer $g\geq2$, let $A_{g}(W)$ be the set of all numbers of the form $\sum\limits_{f\in F}a_{f}g^{f}$ where $F\in \mathcal{F}^{\ast}(W)$ and $1\leq a_{f}\leq g-1$.
 The set $A$ is called an asymptotic basis of order $h$ if $r_h(A,n)\geq 1$ for all sufficiently large integers $n$.
  An asymptotic basis $A$ of order $h$ is minimal if no proper subset of $A$ is an asymptotic basis of order $h$. This means that, for any $a\in A$, the set $E_{a}=hA\setminus h(A\setminus \{a\})$ is infinite.

In 1955, St\"{o}hr \cite{Stohr} first introduced the definition of minimal asymptotic basis. In 1956, H\"{a}rtter \cite{Hartter} gave a nonconstructive proof that there exist uncountably many minimal asymptotic bases of order $h$. In 1988, Nathanson \cite{Nathason1988} proved that if $W_i=\{n\in \mathbb{N}: n\equiv i\pmod h\}$ $(i=0,\ldots,h-1)$, then $\cup_{i=0}^{h-1}A_2(W_i)$ is minimal asymptotic basis of order $h$.
For other related problems, see (\cite{Ling2017}-\cite{Ling2018}).

It is reasonable to consider for any partition $\mathbb{N}=W_0\cup\ldots\cup W_{h-1}$, whether $\cup_{i=0}^{h-1}A_2(W_i)$ is minimal
or not? Nathanson proved this is false even for $h=2$. Moreover, Nathanson posed the following problem(Jia and Nathanson restated this problem again in \cite{Jia1989}).
\begin{problem}\label{prob1} If $\mathbb{N}=W_0\cup\ldots\cup W_{h-1}$ is a partition such that
$w\in W_r$ implies either $w-1\in W_r$ or $w+1\in W_r$, then is
$A=A_2(W_0)\cup\cdots\cup A_2(W_{h-1})$ a minimal asymptotic basis of order $h$?
\end{problem}

In 1989, Jia and Nathanson \cite{Jia1989} obtained the following result.

\noindent{\bf Theorem A} {\it Let $h\geq 2$ and $t=\lceil\log(h+1)/\log 2\rceil$. Partition $\mathbb{N}$ into $h$ pairwise disjoint subsets $W_0,\ldots,
W_{h-1}$ such that each set $W_r$
contains infinitely many intervals of $t$ consecutive integers. Then
$A=A_2(W_0)\cup\cdots\cup A_2(W_{h-1})$ is  a minimal asymptotic basis
of order $h$.}

In 1996, Jia \cite{Jia} generalized Theorem A to $g$-adic representations of integers. In 2011, Chen and Chen \cite{Chen} proved Theorem A under the assumption only required that each set $W_i$ contains one interval of $t$
consecutive integers.

\noindent{\bf Theorem B} {\it Let $h\geq 2$ and $t$ be the least
integer with $t>\log h/\log 2$. Let $\mathbb{N}=W_0\cup\cdots\cup
W_{h-1}$ be a partition such that each set $W_i$ is infinite and
contains $t$ consecutive integers for $i=0,\ldots,h-1$. Then
$A=A_2(W_0)\cup\cdots\cup A_2(W_{h-1})$ is  a minimal asymptotic basis
of order $h$.}

Recently, Yong-Gao Chen and the first author of this paper \cite{ChenTang} proved that the following result:

\noindent{\bf Theorem C} {\it Let $h$ and $t$ be integers with $2\leq t\le \log h/\log 2$.
 Then there exists a partition $\mathbb{N}=W_0\cup\cdots\cup W_{h-1}$
such that each set $W_r$ is a union of infinitely many intervals of at least $t$
consecutive integers and
$$A=A_2(W_0)\cup\cdots\cup A_2(W_{h-1})$$ is not a minimal asymptotic basis of order $h$.
}

It is natural to pose the following $g$-adic version of Problem 1:

\begin{problem}\label{prob1} Let $g\geq 2$ be an integer. If $\mathbb{N}=W_0\cup\cdots\cup W_{h-1}$ is a partition such that
$w\in W_r$ implies either $w-1\in W_r$ or $w+1\in W_r$, then is
$A=A_g(W_0)\cup\cdots\cup A_g(W_{h-1})$ a minimal asymptotic basis of order $h$?
\end{problem}

 Similar to the proof of Theorem B, Ling and Tang remarked in \cite{Ling2015} that Theorem B can be extended to all $g\geq 2$ as following. \vskip 3mm


\noindent{\bf Theorem D} {\it Let $h\geq 2$ and $t$ be the least integer with $t>\max\{1,\frac{\log h}{\log g} \}$, let $\mathbb{N}=W_{0}\cup \cdots \cup W_{h-1}$ be a partition such that each set $W_{i}$ is infinite and contains $t$ consecutive integers for $i=0,1,\ldots,h-1$, then $A=A_{g}(W_{0})\cup \cdots \cup A_{g}(W_{h-1})$ is a minimal asymptotic $g$-adic basis of order $h$.}
\vskip 3mm

 In this paper, we solve the Problem 2.

\begin{theorem}\label{mainthm1} Let $g\geq 2$, $h$ and $t$ be integers with $2\leq t\le \log h/\log g$.
 Then there exists a partition $\mathbb{N}=W_0\cup\cdots\cup W_{h-1}$
such that each set $W_r$ is a union of infinitely many intervals of at least $t$
consecutive integers and
$$A=A_{g}(W_{0})\cup \cdots \cup A_{g}(W_{h-1})$$ is not a minimal asymptotic basis of order $h$.
\end{theorem}
By Theorem D and Theorem 1, we know that the answer to Problem 2 is affirmative for $2\leq h<g^2$ and the answer to Problem 1 is negative for $h\geq g^2$.
For $h>g^t(g-1)$, the following stronger result is proved:

\begin{theorem}\label{mainthm2} Let $g\geq 2$, $h$ and $t$ be integers with $h>g^t(g-1)$.
 Then there exists a partition $\mathbb{N}=W_0\cup\cdots\cup W_{h-1}$
such that each set $W_r$ contains infinitely many intervals of at least $t$
consecutive integers and $n\in hA_g(W_0)$ for all $n\geq h$.
\end{theorem}

\section{Proof of Theorem \ref{mainthm1}}
We need the following lemma:

\begin{lemma}\label{lem1}  (See \cite[Lemma 1]{Lee}.) Let $g\geq2$ be any integer.

$(a)$ If $W_{1}$ and $W_{2}$ are disjoint subsets of $\mathbb{N}$, then $A_{g}(W_{1})\cap A_{g}(W_{2})=\emptyset$.

$(b)$ If $W\subseteq \mathbb{N}$ and $W(x)=\theta x+O(1)$ for some $\theta\in(0,1]$, then there exist positive constants $c_{1}$ and  $c_{2}$ such that $$c_{1}x^{\theta}<A_{g}(W)(x)<c_{2}x^{\theta}$$ for all $x$ sufficiently large.

$(c)$ Let $\mathbb{N}=W_{0}\cup\cdots \cup W_{h-1}$, where $W_{i}\neq \emptyset$ for $i=0,\ldots,h-1$. Then $A=A_{g}(W_{0})\cup\cdots \cup A_{g}(W_{h-1})$ is an asymptotic basis of order $h$.
\end{lemma}

For $a<b$, let
$[a, b]$ denote the set of all integers in the interval $[a, b]$.
Let $\{ m_i \}_{i=1}^\infty$ be a sequence of integers with
$m_1>g^{h+2}$ and $m_{i+1}-m_i>g^{h+2}$ $(i\ge 1)$. Let
$$W_0=[0, m_1] \cup \left( \bigcup_{i=1}^\infty [m_i+t+1, m_{i+1}]  \right) $$
and
$$W_j =\bigcup_{\substack{i=1\\i\equiv j\pmod {h-1}}}^\infty [m_i+1, m_{i}+t],
\quad j=1,2,\dots , h-1.$$ Write
$$A=A_g(W_0)\cup\cdots\cup A_g(W_{h-1}).$$

By Lemma \ref{lem1}, we know that $A$ is an asymptotic basis of order $h$.
It is clear that $g^2\in A_g(W_0)$. Now we prove that $E_{g^2}=hA\setminus
{h(A\setminus \{g^2\})}$ is a finite set. Thus $A$ is not a minimal
asymptotic basis of order $h$.

Let $n>m_2$. We will show that $$n\notin E_{g^2}=hA\setminus
{h(A\setminus \{g^2\})}.$$ This is equivalent to prove that $n\in
h(A\setminus \{g^2\})$.

Let $g$-adic expansion of $n$ be
$$n=\sum_{f\in F_n} a_fg^f, \quad 1\leq a_f\leq g-1.$$
It is clear that $F_n\subseteq \mathbb{N}=W_0\cup W_1\cup \cdots
\cup W_{h-1}$. Divide into the following three cases:

{\bf Case 1:} $F_n\cap W_0= \emptyset $. Then $F_n\subseteq W_1\cup \cdots
\cup W_{h-1}$.

{\bf Subcase 1.1:}  $ |F_n|\ge h$. Then $F_n$ has a partition $$F_n
=L_1\cup L_2\cup \cdots
\cup L_{h},$$ where $L_i\not= \emptyset$ $(1\le i\le h)$ and
for every $L_i$ there exists a $W_j$ $(j\ge 1)$ with $L_i\subseteq
W_j$. Let
$$n_i=\sum_{l\in L_i} a_lg^l, \quad 1\le i\le h.$$
Then $n_i\in A\setminus \{ g^2\} $ and $n=n_1+\cdots +n_{h}$.
Hence $n\in h(A\setminus \{g^2\})$.

{\bf Subcase 1.2:}  $1\le
|F_n|\le h-1$. Write
$$F_n=\{ f_0, \dots , f_{l-1}\}, \quad f_0>\cdots >f_{l-1}.$$
Then $1\leq l\leq h-1$.

Noting that $$n=\begin{cases}\sum\limits_{j=1}^{l-1}a_{f_{j}}g^{f_{j}}+(g-1)\sum\limits_{j=1}^{h-l}g^{f_{0}-j}+g^{f_{0}-(h-l)}, &\text {if }a_{f_{0}}=1,\\
\sum\limits_{j=1}^{l-1}a_{f_{j}}g^{f_{j}}+(a_{f_{0}}-1)g^{f_{0}}+(g-1)\sum\limits_{j=1}^{h-l-1}g^{f_{0}-j}+g^{f_{0}-(h-l-1)}, &\text {if } a_{f_{0}}>1,
\end{cases}$$
moreover, $f_0\ge m_1+1>g^{h+2}>h+2$, thus $f_0-(h-l)>l+2\geq 3$ and $f_0-(h-l-1)>l+3\geq 4$.
Hence $n\in h(A\setminus \{g^2\}).$

{\bf Case 2:} $F_n\cap W_0\not= \emptyset $ and $F_n\cap W_0\not= \{2\} $.

{\bf Subcase 2.1:}  $ |F_n\backslash W_{0}|\ge h-1$. Then $F_n\backslash W_{0}$ has a partition $$F_n\backslash W_{0}
=L_1\cup L_2\cup \cdots
\cup L_{h-1},$$ where $L_i\not= \emptyset$ $(1\le i\le h-1)$ and
for every $L_i$ there exists a $W_j$ $(j\ge 1)$ with $L_i\subseteq
W_j$. Let $L_{0}=F_n\cap W_0$ and
$$n_i=\sum_{l\in L_i} a_lg^l, \quad 0\le i\le h-1.$$
Then $n_i\in A\setminus \{ g^2\} $ and $n=n_0+\cdots +n_{h-1}$.
Hence $n\in h(A\setminus \{g^2\})$.

{\bf Subcase 2.2:}  $1 \leq|F_n\backslash W_{0}|\leq h-2$. Write
$$F_n\backslash W_{0}=\{ f_0, \dots , f_{l-1}\}, \quad f_0>\cdots >f_{l-1}.$$
Then $1\leq l\leq h-2$.
Noting that
$$n=\begin{cases}\sum\limits_{f\in F_n\cap W_0} a_fg^f+\sum\limits_{j=1}^{l-1}a_{f_{j}}g^{f_{j}}+(g-1)\sum\limits_{j=1}^{h-l-1}g^{f_{0}-j}+g^{f_{0}-(h-l-1)}&\text {if }a_{f_{0}}=1,\\
\sum\limits_{f\in F_n\cap W_0} a_fg^f+\sum\limits_{j=1}^{l-1}a_{f_{j}}g^{f_{j}}+(a_{f_{0}}-1)g^{f_{0}}+(g-1)\sum\limits_{j=1}^{h-l-2}g^{f_{0}-j}+g^{f_{0}-(h-l-2)}, &\text {if }a_{f_{0}}>1,
\end{cases}$$
moreover, $f_0\ge m_1+1>g^{h+2}>h+2$, thus $f_0-(h-l-1)>l+3\geq 4$ and $f_0-(h-l-2)>l+4\geq 5$.
Hence $n\in h(A\setminus \{g^2\})$.

{\bf Subcase 2.3:}  $F_n\backslash W_{0}=\emptyset$. That is, $F_n\subseteq W_0$. Write
$$F_n=\{ f_0, \dots , f_{k-1}\}, \quad f_0>\cdots >f_{k-1}.$$
Since $$n>m_2>g^{h+2}>(g-1)(1+g+g^2+\cdots
+g^{h+1}),$$ we have $f_0\ge h+2$.

If $k\geq 3$, then
$$n=\begin{cases}\sum\limits_{j=1}^{k-1}a_{f_{j}}g^{f_{j}}+(g-1)\sum\limits_{j=1}^{h-k}g^{f_{0}-j}+g^{f_{0}-(h-k)}&\text {if }a_{f_{0}}=1,\\
\sum\limits_{j=1}^{k-1}a_{f_{j}}g^{f_{j}}+(a_{f_{0}}-1)g^{f_{0}}+(g-1)\sum\limits_{j=1}^{h-k-1}g^{f_{0}-j}+g^{f_{0}-(h-k-1)}. &\text {if }a_{f_{0}}>1.\end{cases}$$
Hence $n\in h(A\setminus \{g^2\})$.

If $k=2$, then $n=a_{f_{0}}g^{f_{0}}+a_{f_{1}}g^{f_{1}}$.

If $a_{f_{1}}>1$, or $a_{f_{1}}=1, f_{1}\neq 2$, then $$n=\begin{cases}
a_{f_{1}}g^{f_{1}}+(g-1)\sum\limits_{j=1}^{h-2}g^{f_{0}-j}+g^{f_{0}-(h-2)}, & \text{ if }a_{f_{0}}=1\\
a_{f_{1}}g^{f_{1}}+(a_{f_{0}}-1)g^{f_{0}}+(g-1)\sum\limits_{j=1}^{h-3}g^{f_{0}-j}+g^{f_{0}-(h-3)},& \text{ if }a_{f_{0}}>1.
\end{cases}$$
Hence $n\in h(A\setminus \{g^2\})$.

If $a_{f_{1}}=1, f_{1}=2$, then $$n=\displaystyle\begin{cases}(g-1)g+g+(g-1)\sum\limits_{j=1}^{h-3}g^{f_{0}-j}+g^{f_{0}-(h-3)}, & \text{ if }a_{f_{0}}=1\\
(g-1)g+g+(a_{f_{0}}-1)g^{f_{0}}+(g-1)\sum\limits_{j=1}^{h-4}g^{f_{0}-j}+g^{f_{0}-(h-4)}, & \text{ if }a_{f_{0}}>1.
\end{cases}$$
Hence $n\in h(A\setminus \{g^2\})$.

If $k=1$, then
$$n=\displaystyle\begin{cases}(g-1)g^{f_{0}-1}+\cdots+(g-1)g^{f_{0}-(h-1)}+g^{f_{0}-(h-1)}, & \text{ if }a_{f_{0}}=1\\
(a_{f_{0}}-1)g^{f_{0}}+(g-1)\sum\limits_{j=1}^{h-2}g^{f_{0}-j}+g^{f_{0}-(h-2)},& \text{ if }a_{f_{0}}>1.
\end{cases}$$
Hence $n\in  h(A\setminus \{g^2\}).$

{\bf Case 3:} $F_n\cap W_0= \{2\} $.  By $n>m_2$, we have $F_n\setminus
W_0\not=\emptyset$. If $f\in F_n\setminus W_0$, then
$f>m_1>g^{h+2}$. Let
\begin{equation}\label{eq2} n=a_2g^2+\sum_{f\in F_n\setminus \{2\} } a_fg^{f}.\end{equation}

{\bf Subcase 3.1:} $a_{2}>1$ and $|F_n\setminus \{2\}|\geq h-1$.  Then $F_n\backslash \{2\}$ has a partition $$F_n\backslash \{2\}
=L_1\cup L_2\cup \cdots
\cup L_{h-1},$$ where $L_i\not= \emptyset$ $(1\le i\le h-1)$ and
for every $L_i$ there exists a $W_j$ $(j\ge 1)$ with $L_i\subseteq
W_j$. Let $n_{0}=a_2g^2$ and
$$n_i=\sum_{l\in L_i} a_lg^l, \quad 1\le i\le h-1.$$
Then $n_i\in A\setminus \{ g^2\} $ and $n=n_0+\cdots +n_{h-1}$.
Hence $n\in h(A\setminus \{g^2\})$.

{\bf Subcase 3.2:} $a_{2}>1$ and $1\leq |F_n\setminus \{2\}|\leq h-2$. Write
$$F_n\backslash \{2\} =\{ f_0, \dots , f_{l-1}\}, f_0>\cdots >f_{l-1}.$$
Then $1\leq l\leq h-2$.
Noting that
$$n=\displaystyle\begin{cases}a_2g^2+\sum\limits_{j=1}^{l-1}a_{f_{j}}g^{f_{j}}+(g-1)\sum\limits_{j=1}^{h-l-1}g^{f_{0}-j}+g^{f_{0}-(h-l-1)}, & \text{ if }a_{f_{0}}=1\\
a_2g^2 +\sum\limits_{j=1}^{l-1}a_{f_{j}}g^{f_{j}}+(a_{f_{0}}-1)g^{f_{0}}+(g-1)\sum\limits_{j=1}^{h-l-2}g^{f_{0}-j}+g^{f_{0}-(h-l-2)},& \text{ if }a_{f_{0}}>1,
\end{cases}$$
moreover, $f_0\ge m_1+1>g^{h+2}>h+2$, thus $f_0-(h-l-1)>l+3\geq 4$ and $f_0-(h-l-2)>l+4\geq 5$.
Hence $n\in  h(A\setminus \{g^2\}).$

{\bf Subcase 3.3:} $a_{2}=1$ and $|F_n\setminus \{2\}|\geq h-2$.  Then $F_n\backslash \{2\}$ has a partition $$F_n\backslash \{2\}
=L_1\cup L_2\cup \cdots
\cup L_{h-2},$$ where $L_i\not= \emptyset$ $(1\le i\le h-2)$ and
for every $L_i$ there exists a $W_j$ $(j\ge 1)$ with $L_i\subseteq
W_j$. Let $n_{0}=(g-1)g$ and $n_{1}=g$ and
$$n_{i+1}=\sum_{l\in L_i} a_lg^l, \quad 1\le i\le h-2.$$
Then $n_i\in A\setminus \{ g^2\} $ and $n=n_0+\cdots +n_{h-1}$.
Hence $n\in h(A\setminus \{g^2\})$.

{\bf Subcase 3.4:} $a_{2}=1$ and $1\leq |F_n\setminus \{2\}|\leq h-3$. Write
$$F_n\backslash \{2\} =\{ f_0, \dots , f_{l-1}\}, \quad f_0>\cdots >f_{l-1}.$$
Then $1\leq l\leq h-3$. Then $f_0\ge m_1+1>g^{h+2}>h+2$. Hence $f_0-(h-l-2)>l+4\geq 5$ and $f_0-(h-l-3)>l+5\geq 6$.

$$n=\displaystyle\begin{cases}(g-1)g+g+\sum\limits_{j=1}^{l-1}a_{f_{j}}g^{f_{j}}+(g-1)\sum\limits_{j=1}^{h-l-2}g^{f_{0}-j}+g^{f_{0}-(h-l-2)}& \text{ if }a_{f_{0}}=1\\
(g-1)g+g +\sum\limits_{j=1}^{l-1}a_{f_{j}}g^{f_{j}}+(a_{f_{0}}-1)g^{f_{0}}+(g-1)\sum\limits_{j=1}^{h-l-3}g^{f_{0}-j}+g^{f_{0}-(h-l-3)},& \text{ if }a_{f_{0}}>1.
\end{cases}$$
Hence $n\in  h(A\setminus \{g^2\}).$

This completes the proof of Theorem \ref{mainthm1}.

\section{Proof of Theorem \ref{mainthm2}}

Fix an integer $m$ such that $g^m> g^{t+2}h$. Put
$$W_0=\{mk,mk+1,\dots,mk+m-t-1: k\in \mathbb{N}\}$$
and
$$W_{i}=\{mk+m-t,\dots, mk+m-1: k\equiv i\pmod {h-1}\}, \quad 1\le i\le h-1.$$
Let
$$A=A_{g}(W_{0})\cup \cdots \cup A_{g}(W_{h-1}).$$
We will use induction on $n$ to prove that every integer $n\ge h$
is in $hA_g(W_0)$. This implies that $A$ is not a minimal asymptotic
basis of order $h$.

If $n=h$, then by $0\in W_0$ ($k=0$) we have
$n=hg^0\in hA_g(W_0)$.
For $i=1,\ldots, g(g-1)$, we have $i+1\in A_g(W_0)$, so
$n=h+i=(h-1)g^0+(i+1)\in hA_g(W_0)$.

Now we assume that every integer $l$ with  $h\le l<n$ $(n\ge h+g(g-1))$
is in $hA_g(W_0)$. We will prove that $n\in hA_g(W_0)$.

Let $k$ be the integer such that
\begin{equation}\label{21}(g-1)g^{mk}\le n-h<(g-1)g^{m(k+1)}.\end{equation}
Then $k\ge 0$. Let $i$ be the  integer such that
\begin{equation}\label{22}(g-1)g^{mk+i}\le n-h<(g-1)g^{mk+i+1}.\end{equation}
Then $0\le i\le m-1$ and
\begin{equation}\label{23}h\le n-(g-1)g^{mk+i}<(g-1)^2g^{mk+i}+h.\end{equation}
By (\ref{22}) and $n-h\ge g(g-1)$, we have $mk+i\ge 1$. Divide into the following two
cases:

{\bf Case 1:} $0\le i\le m-t-1$.
 By the induction hypothesis,
$n-(g-1)g^{mk+i}\in hA_g(W_0)$. Let
\begin{equation}\label{3.4}n-(g-1)g^{mk+i}=a_1+\cdots +a_h,\quad a_j\in A_g(W_0), j=1,\ldots,h.\end{equation}
If $a_j\ge (g-1)g^{mk+i}$ for all $1\le j\le h$, then by $mk+i\ge 1$ and $h> g^t(g-1)$ we
have
$$n-(g-1)g^{mk+i}=a_1+\cdots +a_h\ge hg^{mk+i}(g-1)>(g-1)^2g^{mk+i}+h,$$
which contradicts with (\ref{23}). So at least one of $a_j$'s is less than $(g-1)g^{mk+i}$.

{\bf Subcase 1.1:} $0\le i\le m-t-2$. Choose an integer $a_j<(g-1)g^{mk+i}$. Write
$$a_j=\sum\limits_{f\in F_0}a_f g^f, \quad 1\leq a_f\leq g-1.$$
Then $f_{\max}:=\max\{f: f\in F_0\}\leq mk+i$.

If $f_{\max}<mk+i$, then $a_j+(g-1)g^{mk+i}\in A_g(W_0)$.

If $f_{\max}=mk+i$, then
 \begin{eqnarray*}a_j+(g-1)g^{mk+i}&=&\sum\limits_{f\in F_0\setminus \{mk+i\}}a_f g^f+a_{mk+i}g^{mk+i}+(g-1)g^{mk+i}\\
 &=&\sum\limits_{f\in F_0\setminus \{mk+i\}}a_f g^f+(a_{mk+i}-1)\cdot g^{mk+i}+g^{mk+i+1},
\end{eqnarray*}
we have $a_j+(g-1)g^{mk+i}\in A_g(W_0).$  Noting that
$$n=(n-(g-1)g^{mk+i})+(g-1)g^{mk+i},$$
by (\ref{3.4}) we have $n\in hA_g(W_0)$.

{\bf Subcase 1.2:} $i=m-t-1$. There exist at least $g-1$ $a_j$'s which are less than $(g-1)g^{mk+m-t-1}$. Otherwise, we have
 \begin{eqnarray*}a_1+\cdots +a_h&\ge&((h-(g-2))(g-1)g^{mk+m-t-1}+(g-2)
\\
&=&(h-(g-1))(g-1)g^{mk+m-t-1}+(g-1)g^{mk+m-t-1}+(g-2)\\
&>&(g-1)^2g^{mk+m-t-1}+h,\end{eqnarray*}
which contradicts with (\ref{23}).
Without loss of generality, we may assume that $$a_{j}<(g-1)g^{mk+m-t-1}, \quad j=1,\ldots,g-1.$$Fix a $j\in\{1,\ldots,g-1\}$, write
$$a_j=\sum\limits_{f\in F_0}a_f g^f, \quad 1\leq a_f\leq g-1.$$
Then $f_{\max}:=\max\{f: f\in F_0\}\leq mk+m-t-1$.

If $f_{\max}<mk+m-t-1$, then $a_j+g^{mk+m-t-1}\in A_g(W_0)$.

If $f_{\max}=mk+m-t-1$, then $1\leq a_{mk+m-t-1}\leq g-2$ and
 \begin{eqnarray*}a_j+g^{mk+m-t-1}&=&\sum\limits_{f\in F_0\setminus \{mk+m-t-1\}}a_f g^f+a_{mk+m-t-1}g^{mk+m-t-1}+g^{mk+m-t-1}\\
 &=&\sum\limits_{f\in F_0\setminus \{mk+m-t-1\}}a_f g^f+(a_{mk+m-t-1}+1)\cdot g^{mk+m-t-1},
\end{eqnarray*}
we have $a_j+g^{mk+m-t-1}\in A_g(W_0).$  Thus
\begin{eqnarray*}n&=&a_1+\cdots +a_h+(g-1)g^{mk+m-t-1}\\
&=&(a_{1}+g^{mk+m-t-1})+\cdots+(a_{g-1}+g^{mk+m-t-1})+a_{g}+\cdots+a_{h}.\end{eqnarray*}
Hence $n\in hA_g(W_0).$

{\bf Case 2:} $m-t\le i\le m-1$. Let $u$ be the  integer such that
\begin{equation}\label{24}h(g-1)g^{mk+u}\le n-h<h(g-1)g^{mk+u+1}.\end{equation} Since $h>g^t$ and $g^m>g^{t+2}h$, we have
\begin{equation}\label{25}g^{mk+i+1}\le g^{mk+m}<hg^{mk+m-t}\end{equation}
and
\begin{equation}\label{26}g^{mk+i}\ge g^{mk+m-t}> hg^{mk+2}.\end{equation}

By (\ref{22}), (\ref{24}) and (\ref{25}), we have
$$h(g-1)g^{mk+u}\le (g-1)g^{mk+i+1}<h(g-1)g^{mk+m-t},$$
thus $u\leq m-t-1$.

By (\ref{22}), (\ref{24}) and (\ref{26}), we have
$$h(g-1)g^{mk+u+1}\ge (g-1)g^{mk+i}>h(g-1)g^{mk+2},$$
thus $u\geq 2$.

So $mk+u\in W_0$ and
$(g-1)g^{mk+u}\in A_g(W_0)$.

{\bf Subcase 2.1:} $u=m-t-1$. Then by (\ref{21}) and $h>g^t(g-1)$, we have
\begin{eqnarray*}n-\sum_{j=2}^{m-t-1}
h(g-1)g^{mk+j}&<&(g-1)g^{mk+m}+h-\sum_{j=2}^{m-t-1} h(g-1)g^{mk+j}\\
&=&(g-1)g^{mk+m}+h-hg^{mk} (g^{m-t}-g^2)\\
&=&g^{mk} ( g^m(g-1-hg^{-t})+g^2h)+h\\
&<& g^{mk} ( g^{t+2}
h(g-1-hg^{-t})+g^2h)+h\\
&=&g^{mk} ( g^2h(g^t(g-1)-h)+g^2h)+h\\
&\le & g^{mk} ( -g^2h+g^2h)+h=h.\end{eqnarray*} Thus
$$n-h<\sum_{j=2}^{m-t-1}
h(g-1)g^{mk+j}.$$ Let
$$\sum_{j=v}^{m-t-1}
h(g-1)g^{mk+j}\le n-h<\sum_{j=v-1}^{m-t-1} h(g-1)g^{mk+j}.$$ Then $3\le v\le
m-t-1$ and
$$0\le n-h-\sum_{j=v}^{m-t-1}
h(g-1)g^{mk+j}<h(g-1)g^{mk+v-1}.$$ Let
$$s(g-1)g^{mk+v-1}\le n-h-\sum_{j=v}^{m-t-1}
h(g-1)g^{mk+j}<(s+1)(g-1)g^{mk+v-1}.$$ Then $0\le s\le h-1$ and
\begin{equation}{\label{3.8}}h\le n-s(g-1)g^{mk+v-1}-\sum_{j=v}^{m-t-1}
h(g-1)g^{mk+j}<(g-1)g^{mk+v-1}+h.\end{equation}
By the induction hypothesis, we may assume that
\begin{equation}{\label{3.9}}n-s(g-1)g^{mk+v-1}-\sum_{j=v}^{m-t-1} h(g-1)g^{mk+j}=b_1+\cdots +b_h, \end{equation}
where $b_j\in A_g(W_0)$, $j=1,\ldots,h$.

By (\ref{3.8}) and (\ref{3.9}), we have \begin{equation}{\label{3.10}}b_{j}<g^{mk+v}, \quad j=1,\ldots,h.\end{equation} Otherwise, we have $$b_1+\cdots +b_h\geq g^{mk+v}+h-1>(g-1)g^{mk+v-1}+h,$$
which contradicts with (\ref{3.8}).
Moreover, we know that there exist at most $g-1$ $b_j$'s which are greater than or equal to $g^{mk+v-1}$. Otherwise, we have
\begin{eqnarray*}b_1+\cdots +b_h&\ge&gg^{mk+v-1}+(h-g)
\\
&=&(g-1)g^{mk+v-1}+g^{mk+v-1}-g+h\\
&>&(g-1)g^{mk+v-1}+h,\end{eqnarray*}
which contradicts with (\ref{3.8}).

Without loss of generality, we may assume that $b_j\geq g^{mk+v-1}$, $j=1,\ldots, l$, where $l\leq g-1$.
Let
$$b_{j}=b_{j}^{\prime}+r_{j}g^{mk+v-1}, \quad j=1,2,\ldots,l$$
be the $g$-adic expansion of $b_j's$ according to the fact that $b_{j}<g^{mk+v}$, $j=1,\ldots,h$. We have $r_{1}+\cdots+r_{l}\leq g-1.$
Otherwise, $$b_1+\cdots+b_h\geq g\cdot g^{mk+v-1}+h-l>(g-1)g^{mk+v-1}+h,$$
which contradicts with (\ref{3.8}).

Noting that $0\le s\le
h-1$, we have
\begin{eqnarray*}n
&=&b_{1}+\cdots+b_{h}+s(g-1)g^{mk+v-1}+\sum_{j=v}^{m-t-1} h(g-1)g^{mk+j}\\
&=&\left(b_{1}^{\prime}+(r_{1}+\cdots+r_{l})g^{mk+v-1}+\sum_{j=v}^{m-t-1} (g-1)g^{mk+j}\right)\\
&+&\left(b_{2}^{\prime}+(g-1)g^{mk+v-1}+\sum_{j=v}^{m-t-1} (g-1)g^{mk+j}\right)\\
&+&\cdots+\left(b_{l}^{\prime}+(g-1)g^{mk+v-1}+\sum_{j=v}^{m-t-1} (g-1)g^{mk+j}\right)\\
&+&\left(b_{l+1}+(g-1)g^{mk+v-1}+\sum_{j=v}^{m-t-1}(g-1)g^{mk+j}\right)\\
&+&\cdots+\left(b_{s+1}+(g-1)g^{mk+v-1}+\sum_{j=v}^{m-t-1} (g-1)g^{mk+j}\right)\end{eqnarray*}
\begin{eqnarray*}&+&\left(b_{s+2}+\sum_{j=v}^{m-t-1} (g-1)g^{mk+j}\right)\\
&+&\cdots+\left(b_{h}+\sum_{j=v}^{m-t-1} (g-1)g^{mk+j}\right).\end{eqnarray*}
Hence $n\in hA_g(W_0)$.

{\bf Subcase 2.2:} $1\le u<m-t-1$. Let
$$s(g-1)g^{mk+u+1} \le n-h<(s+1)(g-1)g^{mk+u+1}.$$
By \eqref{24} and $h>g^t(g-1)$, we have
\begin{equation*} g^{t-1}(g-1)^{2}g^{mk+u+1}<n-h<h(g-1)g^{mk+u+1}.\end{equation*}
Thus $g^{t-1}(g-1)\leq s\leq h-1$.
Let
$$qg^{mk+u}\le n-h-s(g-1)g^{mk+u+1}<(q+1)g^{mk+u}.$$
Noting that
$$0\le n-h-s(g-1)g^{mk+u+1}<g(g-1)g^{mk+u},$$
we have $0\le q\le g(g-1)-1\le h-2$ and
\begin{equation}{\label{3.11}}h\le n-qg^{mk+u}-s(g-1)g^{mk+u+1}<g^{mk+u}+h.\end{equation}
By the induction hypothesis, we may assume that
$$n-qg^{mk+u}-s(g-1)g^{mk+u+1}=c_1+\cdots +c_h, $$
where $c_j\in A_g(W_0)$ $ (1\le j\le h)$. By (3.11) we know that there exists at most one of $c_j$'s is greater than or equal to $g^{mk+u}$.
Hence, we may assume that $$c_j< g^{mk+u}, \quad j=1,\ldots,h$$ or $$c_{1}\geq g^{mk+u}\text{ and } c_j< g^{mk+u}(j=2,\ldots, h).$$
 Noting that $q\leq s$, $g^{t-1}(g-1)\leq s\leq h-1$ and $0\le q\le g(g-1)-1\le h-2$, we have
\begin{eqnarray*} n
&=&c_{1}+\cdots+c_{h}+qg^{mk+u}+s(g-1)g^{mk+u+1}\\
&=&c_{1}+(c_{2}+g^{mk+u}+(g-1)g^{mk+u+1})+\cdots+(c_{q+1}+g^{mk+u}+(g-1)g^{mk+u+1})\\
&+&(c_{q+2}+(g-1)g^{mk+u+1})+\cdots+(c_{s+1}+(g-1)g^{mk+u+1})\\
&+&c_{s+2}+\cdots+c_{h}.
\end{eqnarray*}

Hence $n\in hA_g(W_0)$.

This completes the proof of Theorem \ref{mainthm2}.

\end{document}